\documentclass[twoside]{article}
\usepackage{color,times,amsmath,amsfonts,latexsym,epsfig,graphicx,epsf,colordvi}
\usepackage{url,hyperref}
\usepackage[english]{babel}

\usepackage{amssymb}
\usepackage{comment}
\usepackage{url}
\usepackage{fancyhdr}
\usepackage{mathtools}

\def\tcb{\textcolor[rgb]{0,0,1}}

\newcommand{\CC}{\mathbb {C}}
\newcommand{\RR}{\mathbb {R}}

\newcommand{\bp}{\begin{pmat}}
\newcommand{\ep}{\end{pmat}}

\makeatletter
\def\adots{\mathinner{\mkern1mu\raise\p@\vbox{\kern7\p@\hbox{.}}\mkern2mu\raise4\p@\hbox{.}\mkern2mu\raise7\p@\hbox{.}\mkern1mu}}
\makeatother

\title{ \vspace*{12mm}{ \bf The Construction of High Order Convergent Look-Ahead Finite Difference Formulas for Zhang Neural Networks  \vspace*{-0mm} }}

\author{Frank Uhlig\\
Department of Mathematics and Statistics,\\
Auburn University, AL 36849-5310, USA \ ({\tt uhligfd@auburn.edu})\\[8mm]}

\begin{document}
\date{  }
\maketitle

\vspace*{-2mm}

\pagestyle{fancy}
\fancyhead{}
\fancyhf{} 
\renewcommand{\headrulewidth}{0pt}

\fancyhead[CE]{Frank Uhlig}
\fancyhead[CO]{Convergent look-ahead Difference Formulas} 
\fancyhead[RO]{\thepage}
\fancyhead[LE]{\thepage}
\thispagestyle{empty}

\vspace*{-6mm}
{\normalsize
\noindent 
{\bf Abstract : }
\noindent
Zhang Neural Networks rely on convergent 1-step ahead finite difference formulas of which very few are known. Those which are known have been constructed in ad-hoc ways and suffer from low truncation error orders. This paper develops a constructive method to find convergent look-ahead finite difference schemes of higher truncation error orders. The method consists of seeding the free variables of a linear system comprised of Taylor expansion coefficients followed by a minimization algorithm for the maximal magnitude root of the formula's  characteristic polynomial. This helps us find new convergent 1-step ahead finite difference formulas of any truncation error order. Once a polynomial has been found with roots inside the complex unit circle and no repeated roots on it, the associated look-ahead ZNN discretization formula is convergent and can be used for solving any discretized ZNN based model. Our method  recreates and validates the few known convergent formulas, all of which have truncation error orders at most 4.  It  also creates new convergent 1-step ahead difference formulas with truncation  error orders 5 through 8.\\[2mm]
{\bf Subject Classifications :} \   65Q10, 65L12,  92B20 \\[2mm]
{\bf Key Words :}  finite difference formula, look-ahead difference formula,  Taylor expansion, linear systems, free variable, characteristic polynomial, convergent multistep method,  Zhang neural network,  truncation error order}

\section{\Large Introduction}

Finite difference formulas have a long history of over 200 years in computational mathematics. They came about after the development of Calculus in the late 17th century and were introduced  to estimate the behavior and slope of functions or to approximate areas and volumes. In  the differentiation realm, one of the first such formulas that is still relevant  today is Euler's forward finite difference formula, written here in its symmetric form as
\begin{equation} \label{Eulerappr} \dot{y}_j \approx \dfrac{y_{j+1} - y_{j-1}}{2 \tau} 
\end{equation}
where  $ y_j = y(t_j)$ with $t_j = t_0 + j  \tau$ for a constant step size $\tau$, an initial time instance $t_0$  and any $j \ge 1$. If we solve (\ref{Eulerappr}) for $y_{j+1}$ we obtain the symmetric  1-step ahead finite difference formula of Euler
\begin{equation} \label{EulerFD} y_{j+1}  =  y_{j-1} + 2  \tau  \dot{y}_j \ .
\end{equation}
After assembling the $y_j$ entries on the left of the equation (\ref{EulerFD}) this leads to 
\begin{equation}\label{charpolyE} y_{j+1}  -  y_{j-1} = 2  \tau  \dot{y}_j \ . 
\end{equation}
Next we interpret the left hand side of equation (\ref{charpolyE}) as a polynomial $p$ of smallest degree in a variable $x$ where the subscripts in (\ref{charpolyE}) become the powers of $x$, namely $p(x) = x^2 - 1$. $p$ is called the characteristic polynomial of the difference equation (\ref{charpolyE}). This process is familiar to anyone who has taken a first course in Numerical Analysis and it is described in every elementary textbook on Numerics. We have explained this fundamental process  here in full detail because we need to use it repeatedly in the future.\\[1mm]
 The roots  of the characteristic polynomial $p$ for  Euler's finite difference formula (\ref{EulerFD}) are +1 and --1. They both lie on the periphery of the unit circle in $\CC$ and are distinct. Therefore the 1-step Euler method is convergent. The requirement that all roots of $p$ must lie inside the closed unit disk and no repeated roots may lie on the unit circle is necessary and sufficient for convergence, see e.g.  \cite[ch. 17.6.2]{EMU96}. Convergent finite difference schemes  can be used repeatedly to trace solutions of differential equations in a look-ahead way, albeit for Euler with a low order of accuracy. This was initially done by Bunse-Gerstner et al in \cite{BBMN91} in 1991 for computing time-varying SVDs efficiently via a look-ahead Euler based integrator and subsequently in many other papers.\\[-3mm]

  Zhang Neural Networks were first developed early in the new millennium by Yunong Zhang and others, called  Zeroing Neural Network then, \cite{ZJW2002}. The idea was  taken up by engineers and implemented  in many time-varying applications with well over three hundred articles  published mostly in engineering and applied math journals. Recently there have been impulses from numerical analysts, but ZNN is still largely unknown in numerical circles. ZNN based methods have been used to find time-varying reciprocals, square roots, generalized inverses, pseudoinverses, and  for solving linear, Sylvester, Lyapunov and other equations or inequalities, as well as for matrix eigenvalues and eigenvectors and almost anything matrix related; in optimal design and control and for minimizations and so forth. The list of real-world applications  goes on and on, from robotics to autonomous cars and experimental aircrafts where sensor data arrives at frequencies such as 50 Hz and an objective must be met accurately, 1-step ahead in real-time and for real-world situations. ZNN methods can be easily implemented  in on-chip designs for practical control applications, see e.g. \cite{ZG2015}. 
They currently have many industrial uses and give systems and machines  improved performance. See e.g. 
\cite{LMUZa2018,LMUZb2018,QZY2019,XL2016,ZJW2002,ZKXY2010,ZLYL2012,ZYLHW2017,ZYLUH2018,ZYGZ2011} for a short glimpse of their vast potential with time-varying systems.
Most recently there have been substantial advances on the numerical behavior of continuous-time ZNN methods for matrix problems such as by Lin Xiao et al in \cite{XLLTK} on the stability and robustness in control application, in \cite{XZZLL} on robot applications, and on general time-varying matrix inverse problems  when using time-varying decay constants $\eta(t) > 0$ in \cite{WWS16,WSW18} for example. \\[2mm]
Zhang Neural Networks (ZNN) are  designed for and most efficiently used  to solve time-varying multi-dimensional equations $f(t) = b(t)$ predictively with high accuracy and quickly in real-time. There they use look-ahead finite difference formulas to solve a problem specific error differential equation. All ZNN methods are based  on the error equation $e(t) = f(t) - b(t) \stackrel{!}{=}0$ and the stipulation that $e(t)$ should decay exponentially fast to zero, i.e., 
\begin{equation}\label{eDE}
\dot{e}(t) = -\eta e(t) \text{ for } \eta > 0\ .
\end{equation}
When associated with  a given continuous-time $f(t) = b(t)$ model,  discretized ZNN methods can easily solve the associated discrete time-varying problem $f(t_k) = b(t_k)$ where the time steps $t_k$ are equidistant and the input data is derived from repeated sensor readings. The discretized method predicts the model's solution at  time $t_{k+1}$ from a certain subset of  the previous iterates $f(t_j)$ and $b(t_j)$ with $j \leq k$ and does so shortly after time $t_k$ by using certain derivatives that the error DE (\ref{eDE}) requires. Discrete ZNN methods  must construct the next iterate well before the next time instance $t_{k+1}$ arrives. They succeed  with high accuracy and  a truncation error order of $O(\tau^{m+1})$ if the chosen 1-step ahead discretization formula has truncation error order $m+1$ for the constant sampling gap $\tau$. \\[-3mm]
 
 Section 2 below describes how to set up the mechanics for finding high order convergent look-ahead finite difference formulas via Linear Algebra and elementary Matrix Theory. Section 3 then describes a characteristic polynomial root minimization process that helps us find look-ahead finite difference schemes which satisfy the convergence conditions. Section 4  provides a list of new and high truncation order convergent 1-step ahead finite difference schemes of truncation error orders up to $O(\tau^8)$, as well as open problems. Several of our newly found convergent 1-step ahead finite difference formulas are tested regarding their accuracy and convergence behavior in \cite{FUFoV} on the parameter-varying complex matrix field of values problem.

\section{Construction of general look-ahead Schemes from Seeds via Taylor Polynomials and elementary Linear Algebra}

A recent literature search for  known convergent look-ahead finite difference formulas  found less than half a dozen such formulas, all of which had  truncation error orders less than or equal to 4. This is so despite hundreds of potentially known look-ahead discretization formulas in the literature. Most of them fail the characteristic roots condition and therefore they are unusable for ZNN type methods, such as the look-ahead but unstable formulas (27) to (30) in \cite{LMUZa2018}. With repeated roots on the unit circle, oscillations will eventually set in; with roots outside the unit circle, divergence to infinity will occur whenever corresponding fundamental solutions creep into the current states. With this paper we  more than double the range of available convergent look-ahead methods from truncation error orders  2, 3, and 4,  up to  error orders 5, 6, 7 and 8. The new high error order formulas  speed up convergence and improve accuracy to near machine constant error levels when properly implemented. \\[1mm]
This linear algebraic section  develops the basis for  an algorithm to determine look-ahead finite difference formulas of any truncation order, regardless of convergence or not. The computed look-ahead methods will then be used in the next section to act as initial guesses or seeds for starting a roots minimization process that may find convergent look-ahead schemes or it might not. Note that the previously found convergent methods were all (except for Euler) computed by ad hoc methods with lucky guesses and clever schemes. Here  the process is formalized and developed into a computer code that needs no luck, no sweat and no tears.\\[-3mm]

Let us consider a discrete time-varying state vector $x_j = x(t_j) = x(j \cdot \tau)$ for a constant sampling gap $\tau$ and $j = 0,1,2,...$ and write out $\ell + 1$ explicit Taylor expansions for $x_{j+1}, x_{j-1}, ..., x_{j-\ell}$ around $x_j$ as follows:\\[-5mm]
\begin{eqnarray}
x_{j+1} & =  & x_j ~ +~  \tau \dot{x}_j \ \overbrace { + \ \dfrac{\tau^2}{2!} ~ \ddot{x}_j \ \ ~ + ~ \  \ \dfrac{\tau^3}{3!} ~ \overset{\dots}{x}_j  ~ \ \  ... \ \ \ \ \ \ \ + \  \dfrac{\tau^m}{m!} ~ \overset{m}{\dot{x}}_j \ \ \ \ \ \ \ \ \ } \   ~ + \ O(\tau^{m+1})\\[1mm]
x_{j-1} & =  & x_j ~ - ~ \tau \dot{x}_j \ + \ \dfrac{\tau^2}{2!} ~ \ddot{x}_j \ \ ~ -  ~ \ \ \dfrac{\tau^3}{3!} ~ \overset{\dots}{x}_j  ~ \ \ ...  \ \ \  +  (-1)^m \dfrac{\tau^m}{m!} ~ \overset{m}{\dot{x}}_j \ \ \ + \  O(\tau^{m+1})\\[1mm]
x_{j-2} & =  & x_j - 2\tau \dot{x}_j \ + \dfrac{(2\tau)^2}{2!} ~ \ddot{x}_j -  \dfrac{(2\tau)^3}{3!} ~ \overset{\dots}{x}_j  ~ ...  + (-1)^m \dfrac{(2\tau)^m}{m!} ~ \overset{m}{\dot{x}}_j +O(\tau^{m+1})\\[1mm]
x_{j-3} & =  & x_j - 3\tau \dot{x}_j \ + \dfrac{(3\tau)^2}{2!} ~ \ddot{x}_j -  \dfrac{(3\tau)^3}{3!} ~ \overset{\dots}{x}_j ~  ...  + (-1)^m \dfrac{(3\tau)^m}{m!} ~ \overset{m}{\dot{x}}_j +O(\tau^{m+1})\\
& \vdots& \\[2mm]
x_{j-\ell} & =  & x_j - \ell\tau \dot{x}_j \ \ \underbrace { + \dfrac{(\ell\tau)^2}{2!} ~ \ddot{x}_j -  \dfrac{(\ell\tau)^3}{3!} ~ \overset{\dots}{x}_j ~  ...  \ + (-1)^m \dfrac{(\ell\tau)^m}{m!} ~ \overset{m}{\dot{x}}_j  ~ } + ~ O \ (\tau^{m+1})
\end{eqnarray}
Each right hand side of the Taylor expansion rows or equations above contains $m+2$ terms. The central under- and overbraced $m-1$ 'column terms' on the right hand side of the equal signs each contain a factor of identical powers of $\tau$  and identical  varying order partial derivatives of $x_j$. Our interest  lies only in the remaining 'rational number' factors in  the third through $(m+1)$st  'columns' on the right hand side of equations (5) through (10), i.e., for the moment  we omit the powers of $\tau$ and the  derivatives $\overset{r}{\dot{x}}_j$  for $r = 2, ..., m$ throughout what immediately follows. If we can find a linear combination of the $\ell+1$ equations (5) through (10) that makes the 'braced' $m-1$ number terms in each of these 'columns'  disappear or become 0, we have  found an equation for $x_{j+1}$ in terms of the already known values of $x_j$, $x_{j-1}$, ..., $x_{j-\ell}$, the derivative ${\dot{x}}_j$  with an overall error term of order $O(\tau^{m+1})$. These are the only items that are left once the braced region's linear row combination has become 0. To zero out the braced region,  we now collect the relevant rational numbers factors in the constant matrix $A_{\ell+1,m-1}$ 
\begin{equation}\label{Amatrix}
A = \begin{pmat} \dfrac{1}{2!} & \dfrac{1}{3!} & \dfrac{1}{4!} & \cdots & \dfrac{1}{m!}\\[3mm]
   \dfrac{1}{2!} & -\dfrac{1}{3!} & \dfrac{1}{4!} & \cdots & (-1)^m~ \dfrac{1}{m!}\\[3mm]
   \dfrac{2^2}{2!} & - \dfrac{2^3}{3!} & \dfrac{2^4}{4!} & \cdots & (-1)^m~ \dfrac{2^m}{m!}\\[1.5mm]
     \vdots & \vdots & \vdots & & \\[1mm]
  \dfrac{\ell^2}{2!} & - \dfrac{\ell^3}{3!} & \dfrac{\ell^4}{4!} & \cdots & (-1)^m~ \dfrac{\ell^m}{m!} 
\end{pmat}   \in \RR^{\ell+1,m-1} \ .
\end{equation}
$A$'s entry $a_{u,v}$ in row $u$ and column $v$ is  $\dfrac{(-1)^{v+1} ~ (u-1)^{v+1}}{(v+1)!}$ for $2 \leq u \leq \ell +1$  and $a_{1,v} = \dfrac{1}{(v+1)!}$ for all  $v = 1, ..., m-1$. The complete over- and underbraced summed terms in equations (5) through (10)  has the matrix times vector product form
\begin{equation}\label{ADx}
A_{\ell+1,m-1} \cdot  taudx = \begin{pmat} \dfrac{1}{2!} & \dfrac{1}{3!} & \dfrac{1}{4!} & \cdots & \dfrac{1}{m!}\\[3mm]
   \dfrac{1}{2!} & -\dfrac{1}{3!} & \dfrac{1}{4!} & \cdots & (-1)^m~ \dfrac{1}{m!}\\[3mm]
   \dfrac{2^2}{2!} & - \dfrac{2^3}{3!} & \dfrac{2^4}{4!} & \cdots & (-1)^m~ \dfrac{2^m}{m!}\\[1.5mm]
     \vdots & \vdots & \vdots & & \\[1mm]
  \dfrac{\ell^2}{2!} & - \dfrac{\ell^3}{3!} & \dfrac{\ell^4}{4!} & \cdots & (-1)^m~ \dfrac{\ell^m}{m!} 
\end{pmat}  
\begin{pmat}
\tau^2 \ \ddot{x}_j \\[1mm] \tau^3 \ \overset{\dots}{x}_j \\[1mm] \tau^4 \ \overset{4}{\dot{x}}_j \\[1mm] \vdots \\[1mm]  \tau^{m-1} \ \overset{m-1}{\dot{x}}_j \\[1mm] \tau^m \ \overset{m}{\dot{x}}_j
\end{pmat}
\end{equation}
where the vector $taudx \in \RR^{m-1}$ contains the increasing powers of $\tau $ multiplied by the respective higher derivatives of $x_j$ as entries. If we can find a left kernel row vector $x \in \RR^{\ell+1}$ for $A$ with $x\cdot A = o_{m-1}$, the zero row vector in $\RR^{m-1}$, then $x \cdot A \cdot taudx = 0 \in \RR$ as well. 
If we then form the linear combination of the $\ell+1$ equations in rows (5) through (10) as prescribed by the coefficients of $x$, the linear combination of the sums in  the under- and overbraced columns in   (5) through (10) will vanish and we obtain a single look-ahead formula, involving  only $x_{j+1}, x_j, x_{j-1}, ..., x_{j-\ell}$, and $\dot{x}_j$ with an error term of order  $O(\tau^{m+1})$. A non-zero left kernel vector $x$ for $A_{\ell+1,m-1}$ exists as soon as the number of rows $\ell +1$ of $A$ exceeds the number $m-1$ of columns of $A_{\ell+1,m-1}$. \\[-3mm]

 A left null row vector $x$  for $A$  is a right null column vector for $A^T$ when transposed  and vice versa. 
Column null vectors  can be found  from a reduced row echelon form reduction $R$ of $A^T$ easily by setting the free variables of $R_{m-1,\ell+1}$ equal to a nonzero seed vector $y$ of length $\ell + 1 - (m-1) = \ell -m+2$. In this case -- which will be called the regular case from now on, $R$  is a single row block matrix with the identity matrix $I_{m-1}$ appearing in the first block position  because the coefficient matrix $A$ always has full rank $m-1$ and the linear system is underdetermined. And generally, a dense $m-1$ by $m-1$ matrix $B$ appears in the second block position of $R$, i.e., $R$ has the form 
\begin{equation}\label{Bdef}
R = \begin{pmat} I_{m-1} \ , \ B_{m-1,m-1} \end{pmat}_{m-1,2(m-1)} \ . 
\end{equation}
This  dimensional situation  works very well here. In fact the method works for any matrix $A_{k+s,k}$   and any nonzero seed vector $y \in\RR^s$. Any nonzero seed vector $y \in \RR^s$ spawns a null vector $q \in \RR^{k+s}$ for $A^T$  as  $q =$ {\tt [-R*[zeros(k,1);y];y]} in Matlab notation. We then replace the vector $q$ by $q/q(1)$ in order to arrive at a normalized characteristic polynomial $p$ for the associated convergent look-ahead finite difference equation.\\
Once such a null vector $q\in \RR^{2(m-1)}$ has been computed from a seed vector $y \in \RR^{m-1}$  in the regular $A_{2(m-1),m-1}$ case where $\ell+1 = 2(m-1)$, our algorithm forms the specific linear combination of the set of equations (5) through (10) that $q$ suggests in order to zero out all contributions from the entries in the under- and overbraced third through $m+1$st columns on the right hand side of equations (5) to (10). Then we separate the 1-step ahead state $x_{j+1}$ on the left hand side with  likewise accumulations for the current and earlier  states $x_j, x_{j-1}, ..., x_{j-2m+3}$ and the first time derivative $\dot{x}_j$ according to $q$ on the right hand side of the linearly combined equations. This process  yields  a finite difference multistep formula for  $2m-1$ equidistant state vectors. Its  characteristic polynomial is the normalized polynomial 
 $p =$ {\tt [1;-sum(q);q(2:2(m-1))]}$ \in \RR^{2m-1}$, again in Matlab notation. The polynomial $p$ describes a 1-step-ahead difference equation for $2m\!-\!1$ contiguous instances, i.e., a $(2m\!\!-\!\!1)$-IFD formula in our abbreviation to denote  '$(2m\!\!-\!\!1)$-Instance Difference Formulas' .
\begin{equation}\label{findiffe}
 x_{j+1} + p_2 x_{j} + ... + p_{2m-1} x_{j-2m+3} = c \dot{x}_j \ .
\end{equation}
We  solve (\ref{findiffe}) for $x_{j+1}$ by rearranging  terms and obtain the associated look-ahead rule
\begin{equation}\label{laheadr}
x_{j+1}  = -(p_2 x_{j} + ...  +p_{ 2m -1} x_{j-2m+3}) + c \dot{x}_{\tilde j}    
\end{equation}
where we have incorporated the linear combination of the first derivative $\dot{x}_j$ terms in equations (5) through (10) in the constant $c = p \ .\!* [1;0;-(1:2m-3)']$ (in Matlab notation). Here $p$ (with leading coefficient $p_1$ normalized to 1) is in row vector form  to create the dot product  $c$. And $ \dot{x}_{\tilde j}$ is a backward approximation of the derivative $\dot{x}_{ j} $ at time $t_j$ of sufficiently high truncation error order that can by computed from previous $x_{..}$ state data. To use formula (\ref{laheadr}) in the $A_{2m-2,m-1}$ regular case we clearly need $2m-2$ starting values for the $x_{..}$ terms. And we need $k+s$ starting values $x_{..}$ in the more general $A_{k+s,k}$ situation.
\\[2mm]
Our next task is to find convergent finite difference 1-step ahead formulas of a form such as (\ref{laheadr}). Convergence of such multistep formulas depends on the lay of their characteristic polynomial's roots: these must all lie inside the unit circle in $\CC$ and if there are roots on the unit circle, those must be simple. See e.g. \cite[Ch. 17.6.2 and Definition 17.17, p. 475]{EMU96}. How to achieve convergence behavior here is the subject of our next Section.
\vspace*{-2mm}

\section{A Root Minimization Approach to find convergent look-ahead Difference Schemes of high Error Orders from general look-ahead Schemes}
In Section 2 we have described how to start from a short seed vector and construct a 1-step ahead finite difference scheme for possible use in general ZNN processes. The behavior of the associated  characteristic polynomial determines  convergence or divergence. If the constructed scheme  is not convergent for the chosen seed $y$ then there are roots of its characteristic polynomial that exceed 1 in magnitude or there are repeated roots on the unit circle. In our experience, the former, i.e., roots outside the unit circle, is seemingly always the cause for non-convergence; we have never encountered the latter. Note that the matrix $A_{k+s,k}$  has  very special rational number entries of  integer powers divided by factorials. Therefore the reduced row echelon form $R$ of $A^T$ contains mostly integers or rational numbers with small integer denominators in its second block $B$ as defined in (\ref{Bdef}). \\[1mm]
In this section we propose a minimizing algorithm  for the maximal magnitude characteristic  polynomial root in terms of a given seed vector $y$. For $A_{k+s,k}$ the seed vector space is $\RR^s$ and any seed $y$ therein spawns a unique look-ahead finite difference scheme and an associated characteristic polynomial as was shown in Section 2. However, the set of possible characteristic polynomials itself is not a linear space since sums of such polynomials may or may not be representatives of look-ahead difference schemes. Therefore we can only vary the seed and not the intermediate polynomials in the minimization process and we will have to search indirectly  for a characteristic polynomial with smaller maximal magnitude root in a neighborhood of the starting seed $y$. We implement this minimization process by using  the multidimensional built-in Matlab {\tt fminsearch} minimizer function until we have found a seed with an associated characteristic polynomial that is convergent or there is no convergent such formula from the chosen seed. {\tt fminsearch} uses the Nelder-Mead downhill simplex method  \cite{NM65} that finds local minima for non-linear functions such as ours without using derivatives. It mimics the method of steepest descend and performs a local minimizing search via multiple function evaluations. The main task is to discover generating seed vectors $y$ that  can start the minimizing iteration and find a  characteristic polynomial that is convergent according to the convergent root stipulations. Our seed selection process is currently based on random entry seeds since we know of no better way.\\[1mm]
For the general $A_{k+s,k}$ case and after many different approaches for choosing our random entry seed vectors $ y \in \RR^s$, we decided to start from  seeds $y$ with  normally distributed entries and run {\tt fminsearch} to try and find a local minimum of the maximal magnitude  root of the associated characteristic polynomial. Our minimization algorithm   runs through a double do loop. An outside loop for a number (5 to 20 or 100 ...) of random entry starting seed vectors  and  an inner loop for a number (4 to 7 or 15 ...) of randomized restarts from a previously computed {\tt fminsearch} polynomial that is non-convergent.  We project  its seed onto a  point with newly randomized entries nearby and  use this new seed for a subsequent inner loop run several (2, 5 or 8) times.\\[1mm]
The whole MATLAB code fits onto 80 lines of code, plus 30 lines of comments and two dozen lines of known convergent 1-step ahead discretization formulas of varying truncation error orders between $O(\tau^2)$ to $O(\tau^8)$. Previously  convergent  look-ahead methods were completely unknown for truncation error orders above $O(\tau^4)$. Now we can compute many convergent polynomials of higher error orders quickly.
\vspace*{-2mm}

\section{List of all known and new convergent look-ahead finite Difference\\ Schemes and  some computed Results}

Here is the complete short list of the six known 1-step ahead discretization formulas and their properties ordered by  ascending truncation error order.\\[1mm]
{\bf (A) \  Symmetric Euler Differentiation and Discretization Formula with ZNN truncation error order $O(\tau^2)$} :\\
\hspace*{14mm} $ \dot{y}_j = \dfrac{1}{2\tau} y_{j+1} - \dfrac{1}{2\tau} y_{j-1} + O(\tau)$ or\\[1mm]
\hspace*{14mm} $y_{j+1} = y_{j-1} + O(\tau^2) + ... \text{ problem specific terms from model's  right hand side and } \dot{x}_j $\\[2mm]
\hspace*{6mm} Characteristic Polynomial : $p(x) = x^2 - 1$\\
Formally and according to multistep theory, the Euler formula should lead to a convergent look-ahead discretization formula, but Euler is not convergent as such in practice. Why so is a surprising mystery.\\ [2mm]
{\bf (B) \ 4-IFD Formula from \cite[equations (10), (12)]{LMUZb2018} with ZNN truncation error order $O(\tau^3)$} : \\[1mm]
\hspace*{6mm} $ (10) \ \ \dot{y}_j = \dfrac{1}{\tau} y_{j+1} - \dfrac{3}{2\tau} y_j + \dfrac{1}{\tau} y_{j-1} - \dfrac{1}{2\tau} y_{j-2} + O(\tau^2) $ or\\[1mm]
\hspace*{6mm} $ (12) \ \ y_{j+1} = \dfrac{3}{2} y_j -  y_{j-1} + \dfrac{1}{2} y_{j -2} + O(\tau^3) + ... \text{ problem specific terms from model's rhs and } \dot{x}_j$\\[2mm]
\hspace*{6mm} Characteristic Polynomial : $p(x) = 2x^3 - 3 x^2 + 2x -1$, (not normalized)\\[2mm]
{\bf (C) \ 4-IFD Formula from \cite[equation (11)]{LMUZb2018} with ZNN truncation error order $O(\tau^3)$} : \\[1mm]
\hspace*{6mm} $ (11) \ \ \dot{y}_j = \dfrac{3}{5\tau} y_{j+1} - \dfrac{3}{10\tau} y_j - \dfrac{1}{5\tau} y_{j-1} - \dfrac{1}{10\tau} y_{j-2} + O(\tau^2) $ or\\[1mm]
\hspace*{14mm} $  y_{j+1} = \dfrac{1}{2} y_j + \dfrac{1}{3}  y_{j-1} + \dfrac{1}{6} y_{j -2} + O(\tau^3) + ... \text{ problem specific terms from ...}$\\[2mm]
\hspace*{6mm} Characteristic Polynomial : $p(x) = 6x^3 - 3 x^2 - 2x -1$, (not normalized)\\[2mm]
{\bf (D) \ \ FIFD Formula from \cite[equations (14), (21)]{LMUZa2018} with ZNN truncation error order $O(\tau^3)$} : \\[1mm]
\hspace*{6mm} $ (14) \ \ \dot{y}_j = \dfrac{5}{8\tau} y_{j+1} - \dfrac{3}{8\tau} y_j - \dfrac{1}{8\tau} y_{j-1} - \dfrac{1}{8\tau} y_{j-2}  + O(\tau^2) $ or \\[1mm]
\hspace*{6mm} $ (21) \ \ y_{j+1} = \dfrac{3}{5} y_j + \dfrac{1}{5}  y_{j-1} + \dfrac{1}{5} y_{j -2}  + O(\tau^3) + ... \text{ problem specific terms from ...}$\\[2mm]
\hspace*{6mm} Characteristic Polynomial : $p(x) = 5x^3 - 3 x^2 - x -1$, (not normalized)\\[2mm]
{\bf (E) \ \ 5-IFD Formula from \cite[equations (23), (27)]{LMUZb2018} with ZNN truncation error order $O(\tau^4)$} : \\[1mm]
\hspace*{6mm} $ (23) \ \ \dot{y}_j = \dfrac{4}{9\tau} y_{j+1} + \dfrac{1}{18\tau} y_j - \dfrac{1}{3\tau} y_{j-1} - \dfrac{5}{18\tau} y_{j-2} + \dfrac{1}{9\tau} y_{j-3} + O(\tau^3) $ or \\[1mm]
\hspace*{6mm} $ (27) \ \ y_{j+1} = -\dfrac{1}{8} y_j + \dfrac{3}{4}  y_{j-1} + \dfrac{5}{8} y_{j -2}  - \dfrac{1}{4} y_{j-3} + O(\tau^4) + ... \text{ problem specific terms from ...}$\\[2mm]
\hspace*{6mm} Characteristic Polynomial : $p(x) = 8x^4  +x^3 - 6 x^2 - 5x +2$, (not normalized)\\[2mm]
{\bf (F) \ \ 6N$\tau$CD Formula from \cite[equations (16), (18)]{QZGYL2018} with ZNN truncation error order $O(\tau^4)$} : \\[1mm]
\hspace*{6mm} $ (16) \ \ \dot{y}_j = \dfrac{13}{24\tau} y_{j+1} - \dfrac{1}{4\tau} y_j - \dfrac{1}{12\tau} y_{j-1} - \dfrac{1}{6\tau} y_{j-2} - \dfrac{1}{8\tau} y_{j-3} + \dfrac{1}{12\tau} y_{j-4} + O(\tau^3) $ \ or \\[1mm]
\hspace*{6mm} $ (18) \ \ y_{j+1} = \dfrac{6}{13} y_j  + \dfrac{2}{13}  y_{j-1} + \dfrac{4}{13} y_{j -2}  + \dfrac{3}{13} y_{j-3} - \dfrac{2}{13} y_{j-4} + O(\tau^4) + ... \text{ problem specific terms ...}$ \\[2mm]
\hspace*{6mm} Characteristic Polynomial : $p(x) = 13x^5 -6x^4  -2x^3 - 4 x^2 - 3x +2$, (not normalized)\\[-2mm]

Next we mention a list of new convergent look-ahead discretisation formulas that is detailed within our MATLAB codes. 
We start with some results obtained using our {\tt runconv1step.m}  file \cite{Uconstrcode2018} in Matlab. There are 4 integer inputs to 
{\tt runconv1step.m} : the first input indicates how many outer loop runs are desired   and the second input indicates how many separate repeats with altered seed inputs should be performed in an inner do loop for each outer run. A call of {\tt runconv1step(40,10,k,s)} would thus require 40 outer loop runs with 10 separate seeds each (400 runs of {\tt fminsearch} in total). This call tries to find convergent  polynomials of truncation error orders $k+2$ , i.e., polynomials of degree $k+s$ with roots properly inside  the closed unit disk in $\CC$. Here $k$ can be any integer less than 6 and $s$ should be at least equal to $k$ so that the rational entry matrix $A_{k+s,k}$ has a nontrivial left nullspace. Here $s$ is the number of real entries in the seed vector $y$ and our   highest tested-for truncation error order $O(\tau^{k+2})$ is  $k+2 = 6+2 =8$ for $k = 6$. We see no real need to go beyond $k = 6$ since a truncation error order of $(50 Hz)^8 = 0.02^8 \approx 2.56 \cdot 10^{-14}$ seems close enough to the machine constant to wonder about further improvements.\\[-3mm]

The MATLAB output for validating the known convergent 5-IFD formula (E) above from \cite{LMUZb2018} which is listed in the examples list  in our code {\tt runconv1step.m} \ \cite{Uconstrcode2018} is as follows with $y$ denoting the seed vector of length $s$ that was used.

\begin{verbatim}
format short, tic, TPOLY =  runconv1step(1,1,2,2), toc

Truncation_error_order =
     4
y =
    -5     2
TPOLY = 1.0000  0.1250 -0.7500 -0.6250  0.2500  0  -0.0000  0.9025  2.2500

Elapsed time is 0.008423 seconds.
\end{verbatim}

\noindent
The first $k+s+1$ entries in {\tt TPOLY} are the  computed coefficients of the normalized convergent polynomial $p(x) = x^4  +0.125x^3 - 0.75 x^2 - 0.625x + 0.25$ of degree $k+s = 4$ in decreasing exponent order. These are followed by a data separating zero, the deviation of the maximal magnitude root of $p$ from 1 (which is nearly zero), the magnitude of the second largest magnitude root of $p$ and finally the coefficient of $\tau$ that is to be used inside the ZNN discretization. The reader can sample any of the ZNN papers in the bibliography to learn how to set up and implement a  discrete ZNN method from any known convergent polynomial.\\
Next a similar example  for a convergent truncation error order 5  formula from the seed {\tt y = [a,110,-40]} with variable constant $1 \leq a \leq 2.5$, also given in  our  formula list inside  {\tt runconv1step.m}. We have set   $a = 1$ to obtain the output below.\\[-4mm]
\begin{verbatim}
format rat, tic, TPOLY =  runconv1step(1,1,3,3), toc
 
Truncation_error_order =
       5       
y =
       1            110            -40       
TPOLY = 
 1 80/237 -182/237 -206/237 1/237 110/237 -40/237 0 0.0000 446/465 196/79 
   
Elapsed time is 0.007906 seconds.
\end{verbatim}
Note that the resulting normalized characteristic polynomial has again only rational coefficients. This occurrence is rather rare. 
If we run {\tt runconv1step(20,10,5,5)} with $20 \cdot 10 = 200$  individual searches we typically capture 4 to 6 convergent polynomials with discrete ZNN model truncation error orders $O(\tau^7)$. When doing the same with $k = 4$ and $s = 4$ in 200 searches, our code discovers  around 40 convergent polynomials with associated truncation error order $O(\tau^6)$ in ZNN applications. It seems  advantageous to increase $s$ beyond $s = k$ to find  more higher truncation error order polynomials quickly. For example  {\tt runconv1step(20,10,5,6)} with $s = k+1 = 6$ exhibits around a dozen good polynomials and {\tt runconv1step(20,10,5,7)} with $s$ set to $k+2 = 7$ around 18. These success numbers are influenced by the random nature of our seeds $y$ of any fixed length $s\geq k$. Note that if the call of  {\tt runconv1step(runs,jend,k,s)}  generates feasible polynomials, these will all be of degree $k+s$. Therefore an implementation inside any discrete ZNN model must compute $k+s$ starting values before  an associated 1-step forward discrete ZNN iteration can be run.\\[1mm]
We have found  very few high order convergent polynomials with integer coefficients. The usual output {\tt TPOLY} from {\tt runconv1step(...,...,k,s)}  comes  in 16 digit exponent 10 Matlab notation. The look-ahead convergent finite difference formulas of high truncation error orders $O(\tau^5)$ to $O(\tau^8)$  at the bottom of our  our list in 
{\tt runconv1step(...,...,k,s)}  can be generated from the given seed information there. We have not tested the quality of the newly computed convergent polynomials and whether there are any  noticeable differences in real-world problems between them. For example  we do not know whether those polynomials with well separated largest and second largest magnitude roots perform better than those with multiple very near 1 magnitude  roots. To study this issue,  we have included the magnitude of the second largest root in the last but one column of {\tt TPOLY}.\\[1mm]
This paper describes how we can now create high error order convergent look-ahead finite difference schemes for any $k$ and $s \geq k$. Our newly found  \verb+k_s+ designated polynomials \verb+3_3+, \verb+4_5+, \verb+5_5+ and \verb+5_6+ of truncation orders 5 through 7 were tested recently on the matrix field of values problem and compared in \cite{FUFoV}.\\[-2mm]

We conclude with a more theoretical aspect and an open question encountered with finding convergent  polynomials for 1-step ahead ZNN processes.\\[-2mm]

\noindent
 {\bf Remark} and {\bf Open Question :} \\[2mm]
\hspace*{5.5mm} \begin{minipage}{154mm}{ Every polynomial that we have constructed from any seed vector $y \in \RR^s$  by our method has had at least one root on the unit circle (within $10^{-15}$ numerical accuracy).   This is so even for  non-convergent polynomials with some roots outside the unit disk.  Is it true  in general that all  such Taylor expansion matrix $A$ based polynomials have at least one root  on the unit circle in $\CC$?  Or are there some such convergent polynomials whose roots all lie inside the open unit disk.} 
\end{minipage}

\section{Assessment and Outlook}
Our method to construct convergent look-ahead finite difference equations
hinges on two separate ideas and subsequently two constructive steps. The first problem, part one, is essentially linear. We want to eliminate the second to $m$th derivatives in the Taylor expansion equations (5) to (10) around $x_j$ and find a  certain  linear combination of the $\ell +1$ equations that can give us a candidate finite difference scheme. This linear first branch of our task starts from a short seed vector $y$. It completes the seed  of length $s$ to a full  difference equation of length $k+s$ and the associated characteristic polynomial coefficients. The resulting difference equation  relates the 1-step ahead $x_{j+1}$ value  to earlier $x_k$ values with $k \leq j$ and to the derivative at $x_j$ with an error term of order $O(\tau^{m+1})$. The second problem, part two,  is highly nonlinear. It tries to select those finite difference equations whose characteristic polynomials satisfy  certain root conditions for convergence. To solve the second, the nonlinear polynomial roots problem, we have chosen a multidimensional minimization function that starts from a candidate finite  difference scheme an its associated set of coefficients and varies the original seed with an eye on minimizing the largest magnitude root of the associated varying characteristic polynomial. It is the nature of the very first seed from the linear part one that eventually determines the convergence qualities of the method  at the end of each minimizing part two  process. Only then do we know.\\[1mm]
In most of our trial test runs from a seed to a possibly convergent finite difference scheme, the original maximal magnitude characteristic polynomial root does not dip down to 1 or below. Instead the maximal magnitude root usually settles around  1.008.., 1.01.. or larger and does no longer budge, indicating that the original seed vector is not allowing our algorithm to find a usable difference scheme and that we must end this run. Clearly not all  seeds can solve our problem, their surrounding  maximal root 'valleys' simply may not dip low enough to be of use to us. But often enough, our randomized seed minimization algorithm leads to success in finding convergent look-ahead finite difference formulas with truncation error orders up to 8 where none were previously known.\\[2mm]

 \vspace*{5mm}

\hspace*{30mm}{[} \ .../latex/constructconv1stepahead.tex \ ]  \hfill    \today



\end{document}